\documentclass[twoside,reqno,times]{amsart}
\usepackage[active]{srcltx}
\usepackage{latexsym,rotate,eucal,cite}
\usepackage{amsmath,amsthm,amssymb,amsxtra}

\newtheorem{dl}{Theorem}[section]
\newtheorem{yl}[dl]{Lemma}

\theoremstyle{definition}
\newtheorem{dy}[dl]{Definition}
\newtheorem{tl}[dl]{Corollary}
\newtheorem{lz}[dl]{Example}

\newtheorem{xinzhi}[dl]{Proposition}
\newtheorem{remark}[dl]{Remark}

\newcommand{\poq}[2]{(#1;q)_{#2}}

\def\qed{\hfill \rule{4pt}{7pt}}
\def\pf{\noindent {\it Proof.} }

\numberwithin{equation}{section}

\begin{document} 
\title{A general $q$-expansion formula based on matrix inversions and its applications}
       \dedicatory{    \textsc{Jin Wang~\dag}\\[1mm]
            Department of Mathematics, Soochow University\\
            Suzhou 215006, P.~R.~China\\
Email:~\emph{jinwang2016@yahoo.com}
\thanks{\dag~Corresponding author. This work was supported by NSFC (Grant~No. 11471237)}}


\begin{abstract}
In this paper, by use of matrix inversions, we establish a general $q$-expansion formula of arbitrary formal power series $F(z)$ with respect to the base
 $$
\left\{z^n\frac{\poq{az}{n}}{\poq{bz}{n}}\bigg|n=0,1,2\cdots\right\}.$$
 Some concrete expansion formulas and their applications to $q$-series identities are presented, including Carlitz's $q$-expansion formula, a new partial theta function identity and a coefficient identity of Ramanujan's ${}_1\psi_1$ summation formula as  special cases.
\end{abstract}

\maketitle\thispagestyle{empty}
\markboth{ J. Wang}{A general $q$-expansion formula based on matrix inversions and its applications}

\section{Introduction}
Throughout the present paper, we adopt the standard notation and
 terminology for  $q$-series from the book  \cite{10}.  As customary, the $q$-shifted factorials of complex variable $z$
 with the base $q: |q|<1$ are  given by
\begin{eqnarray}
(z;q)_\infty
:=\prod_{n=0}^{\infty}(1-zq^n)\qquad\mbox{and}\quad (z;q)_n:=\frac{(z;q)_\infty}{(zq^n;q)_\infty}\label{guiding}
\end{eqnarray}
for all integers $n$. For integer $m\geq 1$, we use
 the multi-parameter compact notation
\[(a_1,a_2,\ldots,a_m;q)_n:=(a_1;q)_\infty (a_2;q)_n\ldots (a_m;q)_n.\]
Also,  the ${}_{r+1}\phi_r$ series with the
base $q$ and the argument $z$ is defined to be
\begin{align*}
{}_{r+1}\phi_r\bigg[\genfrac{}{}{0pt}{}{a_1,a_2,\cdots,a_{r+1}}{b_1,b_2,\cdots,b_{r}};q,z\bigg]&:=\sum _{n=0} ^{\infty}
\frac{\poq{a_1,a_2,\ldots,a_{r+1}}{n}}{\poq{q,b_1,b_2,\ldots,b_r}{n}}z^{n}.
\end{align*}
For any   $f(z)=\sum_{n\geq 0}a_nz^n\in \mathbb{C}[[z]]$, $\mathbb{C}[[z]]$ denotes  the ring  of  formal power series   in variable $z$,  we shall employ the coefficient functional
$$\boldsymbol\lbrack z^n\boldsymbol\rbrack \{f(z)\}:=a_n\,\,\mbox{and}\,\, a_0=f(0).$$
We also follow the summation convention that for any integers $m$ and $n$,  $$\sum_{k=m}^{n}a_k=-\sum_{k=n+1}^{m-1}a_k.$$

 In their paper \cite{ono}, G.H. Coogan and K. Ono presented the following identity which leads to the generating functions  for values of
certain expressions of Hurwitz zeta functions at non-positive integers.
\begin{yl}[cf.\mbox{\cite[Proposition 1.1]{ono}}]\label{ma-id11} For $|z|<1$, it holds
\begin{align}
\sum_{n=0}^\infty\,z^n\frac{\poq{z}{n}}{\poq{-z}{n+1}}=
\sum_{n=0}^{\infty}(-1)^nz^{2n}q^{n^2}.\label{id11}
\end{align}
\end{yl}
The appearance of \eqref{id11} reminds us of the famous Rogers--Fine identity \cite[Eq. (17.6.12)]{dlmf}.
\begin{yl}For $|z|<1$, it holds
\begin{align}
(1-z)\sum_{n=0}^\infty\,z^n\frac{\poq{aq}{n}}{\poq{bq}{n}}=
\sum_{n=0}^{\infty}(1-azq^{2n+1})(bz)^nq^{n^2}\frac{\poq{aq,azq/b}{n}}{\poq{bq,zq}{n}}.\label{id117}
\end{align}
\end{yl}
In fact, Identity \eqref{id11} can be easily deduced from
\eqref{id117} via setting $aq=z=-b$.
Moreover, by setting $a=z=-b$ in \eqref{id117}, we obtain another similar identity.
\begin{yl}\label{ma-id1177}For $|z|<1$, it holds
\begin{align}
\sum_{n=0}^\infty\,z^n\frac{\poq{z}{n+1}}{\poq{-zq}{n}}=1+
2\sum_{n=1}^{\infty}(-1)^nz^{2n}q^{n^2}.\label{id1177}
\end{align}
\end{yl}
It is these identities, once treated  as formal power series in $z$, that make us be aware  of  investigating in a full generality the problem of representations of formal power series in terms of the sequences
$$
\left\{z^n\frac{\poq{az}{n}}{\poq{bz}{n}}\bigg|n=0,1,2\cdots\right\},$$ which is just a  base of the ring $\mathbb{C}[[z]]$.  This fact asserts  that for  any  $F(z)\in \mathbb{C}[[z]]$, there exists  the series expansion
\begin{align}
F(z)=\sum_{n=0}^\infty\,c_nz^n\frac{\poq{az}{n}}{\poq{bz}{n}},\label{import1-0}
\end{align}
 where the coefficients $c_n$ must be independent of $z$ but may depend on the parameters $a$ and $b$. In this respect,  particularly
noteworthy is  that in \cite{xrma}, X.R. Ma established a (formal) generalized Lagrange--B\"{u}rmann inversion formula.
We record it for direct reference.
\begin{yl}[cf.\mbox{\cite[Theorem 2.1]{xrma}}]\label{analogue}
 Let  $\{\phi_n(z)\}_{n=0}^{\infty}$ be arbitrary sequence of formal power series with $\phi_n(0)=1$. Then for any $F(z)\in \mathbb{C}[[z]]$, we have
 \begin{subequations}\label{expan-one}
\begin{equation}
F(z)=\sum_{n=0}^{\infty}\frac{c_nz^{n}}{\prod_{i=0}^{n}\phi_i(z)},\label{expan15a}
\end{equation}
where the coefficients
\begin{eqnarray}
c_n=\sum_{k=0}^n\boldsymbol\lbrack z^{k}\boldsymbol\rbrack\{F(z)\}\sum_{\stackrel{i\leq j_i}{k=j_{0}\leq j_1\leq j_2\leq\cdots\leq j_{n}
\leq j_{n+1}=n}}\prod_{i=0}^n\boldsymbol\lbrack z^{j_{i+1}-j_{i}}\boldsymbol\rbrack\{\phi_i(z)\}.\label{expan-one-ii}
\end{eqnarray}
\end{subequations}
\end{yl}
For further information on Lemma \ref{analogue}, we refer the reader to \cite{xrma}. As for the classical  Lagrange--B\"{u}rmann inversion formula  the reader might consult the book  \cite[p. 629]{andrews4}  by G.E. Andrews, R. Askey, and R. Roy.  For its various $q$-analogues, we refer the reader to the paper  \cite{andrews} by G.E. Andrews, \cite{carliz} by L. Carlitz, \cite{gessel-stan} by I. Gessel and D. Stanton, and \cite{kratt} by Ch. Krattenthaler,  especially to the good survey  \cite{11} of D. Stanton for a more comprehensive information.

A simple expression of the coefficients $c_n$  seems unlikely  under the case \eqref{expan15a}. Without doubt,  such an explicit  formula   is the key step to successful use  of this expansion formula.
But in contrast,  as far as \eqref{import1-0} is concerned,  we are able to establish  the following explicit expression of $c_n$ via the use of matrix inversions (see Definition 2.1 below). It is just the theme of  the present paper.
\begin{dl}\label{analogue-two}
For  $F(z)\in \mathbb{C}[[z]]$, there exists the series expansion
 \begin{subequations}\label{expan-two}
\begin{equation}
F(z)=\sum_{n=0}^{\infty}c_nz^{n}\frac{\poq{az}{n}}{\poq{bz}{n}}\label{expan-two-1}
\end{equation}
with the coefficients
\begin{align}
c_n=\boldsymbol\lbrack z^{n}\boldsymbol\rbrack\bigg\{  F(z)\frac{\poq{bz}{n-1}}{\poq{az}{n}}\bigg\}-a\sum_{k=0}^{n-1}B_{n-k,1}(a,b)q^{(n-k)k}\boldsymbol\lbrack z^{k}\boldsymbol\rbrack\bigg\{F(z)\frac{\poq{bz}{k}}{\poq{az}{k+1}}\bigg\},\label{coeformula}
\end{align}
where $B_{n,1}(a,b)$ are given by
\begin{align}
z=\sum_{n=1}^\infty B_{n,1}(a,b)z^n\frac{\poq{az}{n}}{\poq{bz}{n}}.\label{180}
\end{align}
\end{subequations}
\end{dl}
As a direct application of this expansion formula, we further set up a general transformation concerning  the Rogers--Fine identity \eqref{id117}.
\begin{dl}\label{tlidentity} For  $G(z)=\sum_{n=0}^{\infty}t_nz^n\in \mathbb{C}[[z]]$, it always holds
\begin{subequations}
\begin{align}
\frac{\poq{az}{\infty}}{\poq{bz}{\infty}}G(z)&=
\sum_{n=0}^{\infty}\frac{\poq{aq/b,az}{n}}{\poq{q,bz}{n}}(bz)^nq^{n(n-1)}\nonumber\\
&\qquad\times\big(\widetilde{G}(zq^n;a,b)-
azq^{2n}\widetilde{G}(zq^{n+1};a/q,b/q)\big),
\label{expan-two-lzlzlz}
\end{align}
where
\begin{align}
\widetilde{G}(z;a,b)&:=\sum_{n=0}^\infty t_nz^n\frac{\poq{az}{n}}{\poq{bz}{n}}.\label{dy}
\end{align}
\end{subequations}
\end{dl}

The rest of this paper is organized as follows. In Section 2, we shall  prove Theorem \ref{analogue-two}. For this purpose, a series of preliminary results will be established. Section 3 is devoted to the proof of Theorem\ref{tlidentity}. Some applications of these two theorems to $q$-series are further  discussed. Among these applications, there is a new partial theta function identity and a coefficient identity of Ramanujan's ${}_1\psi_1$ summation formula.

\section{The proof of Theorem \ref{analogue-two}}
In this section, we proceed to show  Theorem \ref{analogue-two} which amounts to finding  the coefficients $c_n$. For this purpose, we need the concept of matrix inversions and a series of  preliminary  lemmas.
\begin{dy}(cf.\cite[Chapters 2 and 3]{riodan} or
\cite[Definition 3.1.1]{egobook})\label{ddd000} A pair of infinite lower-triangular matrices  $A=(A_{n,k})$ and $B=(B_{n,k})$  is  said to be inverses to each other if and only if for any integers $n, k\geq 0,$
\begin{align}
\sum_{i=k}^n\,A_{n,i}B_{i,k}=\sum_{i=k}^n\,B_{n,i}A_{i,k}
=\left\{
\begin{array}{ll}
0, & n\neq k; \\
 1, & n=k,
\end{array}
\right.
\end{align}
where $A_{n,k}=B_{n,k}=0$ if $n<k.$
As usual, we also say that $A$ and $B$ are invertible and write $A^{-1}$ for $B$.
\end{dy}
Consider  now a particular matrix $A=(A_{n,k})$ with the entries $A_{n,k}$ given by
\begin{align}
z^k\frac{\poq{az}{k}}{\poq{bz}{k}}=\sum_{n=k}^\infty\,A_{n,k}z^n.\label{eq1-13}
\end{align}
It is easy to check that $A=(A_{n,k})$ is invertible. In what follows, let us assume  its inverse
$$A^{-1}=(B_{n,k}(a,b)).$$ As such, we see that \eqref{eq1-13}
is equivalent to
\begin{align}
z^k=\sum_{n=k}^\infty B_{n,k}(a,b)z^n\frac{\poq{az}{n}}{\poq{bz}{n}}.\label{18}
\end{align}

Next, we shall focus on two kinds of generating functions of the entries  $B_{n,k}(a,b)$ of the  matrix $A^{-1}$.
\begin{yl}\label{ddd} Let $B_{n,k}(a,b)$ be the same as above.
 Then we have
\begin{align}B_{n,k+1}(a,b)+(b-a)\sum_{i=k+2}^{n}b^{i-k-2}  B_{n,i}(a,b)=q^{n-k-1}B_{n-1,k}(a,b).\label{impor-1-added}
\end{align}
\end{yl}
\pf
At first, by replacing $z$ by $zq$ in  \eqref{18}, we have
\begin{align*}
(zq)^k=\sum_{n=k}^\infty B_{n,k}(a,b)(zq)^n\frac{\poq{azq}{n}}{\poq{bzq}{n}}.
\end{align*}
Multiplying both sides with $z(1-az)/(1-bz)$ and shifting $n$ to $n-1$, we obtain
\begin{align}
q^k z^{k+1}\bigg(1+(b-a)z\sum_{i= 0}^{\infty}b^{i}z^{i}\bigg)=\sum_{n=k+1}^\infty B_{n-1,k}(a,b)q^{n-1}z^n\frac{\poq{az}{n}}{\poq{bz}{n}}.\label{1888}
\end{align}
Viewing \eqref{1888} from the definition of \eqref{18}, we find that \eqref{1888} can be recast as
\begin{align}
&q^k\sum_{n=k+1}^\infty B_{n,k+1}(a,b)z^n\frac{\poq{az}{n}}{\poq{bz}{n}}\nonumber\\
&+(b-a)q^k\sum_{i= 0}^{\infty}b^{i}\sum_{n=k+i+2}^\infty B_{n,k+i+2}(a,b)z^n\frac{\poq{az}{n}}{\poq{bz}{n}}
\nonumber\\
&=\sum_{n=k+1}^\infty B_{n-1,k}(a,b)q^{n-1}z^n\frac{\poq{az}{n}}{\poq{bz}{n}}.\label{18888}
\end{align}
Thus, by the uniqueness of  the coefficients under the base $\{z^n\poq{az}{n}/\poq{bz}{n}\}_{n=0}^\infty$, it holds
\begin{align*}
q^k B_{n,k+1}(a,b)+(b-a)q^k\sum_{i=0}^{n-k-2}b^{i}  B_{n,k+i+2}(a,b)=q^{n-1}B_{n-1,k}(a,b).
\end{align*}
After slight simplification, we obtain
\begin{align*}
B_{n,k+1}(a,b)+(b-a)\sum_{i=k+2}^{n}b^{i-k-2} B_{n,i}(a,b)=q^{n-k-1} B_{n-1,k}(a,b).
\end{align*}
Hence \eqref{impor-1-added} follows.
\qed

By use of Lemma \ref{ddd}, it is easy to set up a bivariate generating function of  $\{B_{n,k}(a,b)\}_{n\geq k\geq 0}$.
\begin{yl} \label{bgf} Let $B_{n,k}(a,b)$ be defined by \eqref{18}. Then  we have
\begin{align}
G(y,z)=\sum_{n=0}^\infty\frac{\poq{b/y}{n}}{\poq{a/y}{n}}(yz)^n+\sum_{n=0}^\infty (bzq^n-aG_1(zq^n))\frac{\poq{b/y}{n}}{\poq{a/y}{n+1}}(yz)^n,\label{impor-1}
\end{align}
where
\begin{align}
G(y,z)&:=\sum_{k=0}^\infty G_k(z)y^k,\label{impor-2}\\
G_k(z)&:=\sum_{n=k}^\infty B_{n,k}(a,b)z^n.\label{impor-3}
\end{align}
\end{yl}
\pf It suffices to multiply both sides of \eqref{impor-1-added} with $b$. Then
 we get
\begin{align}bB_{n,k+1}(a,b)+(b-a)\sum_{i=k+2}^{n}b^{i-k-1}  B_{n,i}(a,b)= bq^{n-k-1}B_{n-1,k}(a,b).\label{sansan-5}
\end{align}
Shifting $k$ to $k-1$  in \eqref{impor-1-added} gives rise to
\begin{align}
B_{n,k}(a,b)+(b-a)\sum_{i=k+1}^{n}b^{i-k-1}  B_{n,i}(a,b)= q^{n-k}B_{n-1,k-1}(a,b).\label{sansan-6}
\end{align}
By abstracting \eqref{sansan-5} from \eqref{sansan-6}, we come up with
\begin{align}
B_{n,k}(a,b)-aB_{n,k+1}(a,b)=q^{n-k}B_{n-1,k-1}(a,b)-bq^{n-k-1}B_{n-1,k}(a,b).\label{sansan-007}
\end{align}
After multiplying \eqref{sansan-007} by $z^{n}$ and summing over $n$ for $n\geq k$, then we have
\begin{align*}
\sum_{n=k}^\infty B_{n,k}(a,b)z^{n}&-a\sum_{n=k}^\infty B_{n,k+1}(a,b)z^{n}\\
&=zq^{1-k}\sum_{n=k}^\infty B_{n-1,k-1}(a,b)(qz)^{n-1}-bzq^{-k}\sum_{n=k}^\infty B_{n-1,k}(a,b)(qz)^{n-1}.
\end{align*}
In terms of $G_k(z)$ defined by \eqref{impor-3}, this relation can be expressed  as
\begin{align}
G_k(z)-aG_{k+1}(z)=zq^{1-k}G_{k-1}(qz)-bzq^{-k}G_k(qz).\label{sansan-77777}
\end{align}
Then, by multiplying  $y^{k+1}$ and then summing over $k$ for $k\geq 1$ on both sides of \eqref{sansan-77777}, we further obtain
\begin{align*}
(y-a)G(y,z)+(a-y)G_0(z)+ayG_1(z)=yz(y-b)G(y/q,qz)+byzG_0(qz),
\end{align*}
where $G(y,z)$ is given by \eqref{impor-2}. Observe that $G_0(z)=1$.
Then
\begin{align*}
(y-a)G(y,z)-yz(y-b)G(y/q,qz)=y-a+byz-ayG_1(z),
\end{align*}
namely,
\begin{align}
G(y,z)-yz\frac{1-b/y}{1-a/y}G(y/q,zq)=d(y,z),\label{123}
\end{align}
where, for clarity, we define
\[d(y,z):=1+\frac{y}{y-a}(bz-aG_1(z)).\]
Setting $(y,z)\to(y/q^n,zq^n)$ in \eqref{123}, we obtain its  equivalent version as below
\begin{align}
X(n)-yz\frac{1-bq^n/y}{1-aq^n/y}X(n+1)=d(y/q^n,zq^n),\label{123-123}
\end{align}
where
\[X(n):=G(y/q^n,zq^n).\]
Iterating \eqref{123-123} $m$ times, we find
\begin{align*}
X(n)&=d(y/q^n,zq^n)+yz\frac{1-bq^n/y}{1-aq^n/y}X(n+1)\\
&=d(y/q^n,zq^n)+yz\frac{1-bq^n/y}{1-aq^n/y}d(y/q^{n+1},zq^{n+1})\\
&\qquad+(yz)^2\frac{(1-bq^n/y)(1-bq^{n+1}/y)}{(1-aq^n/y)(1-aq^{n+1}/y)}X(n+2)\\
&=\cdots\\
&=\sum_{k=0}^{m-1} d(y/q^{n+k},zq^{n+k})\frac{\poq{bq^n/y}{k}}{\poq{aq^n/y}{k}}(yz)^k+\frac{\poq{bq^n/y}{m}}{\poq{aq^n/y}{m}}(yz)^mX(n+m).
\end{align*}
Regarding the solution of this recurrence relation, we may guess and then show by induction on $m$ (set $n=0$) that
\begin{align*}
G(y,z)&=\sum_{k=0}^\infty d(y/q^k,zq^k)\frac{\poq{b/y}{k}}{\poq{a/y}{k}}(yz)^k\\
&=\sum_{k=0}^\infty\frac{\poq{b/y}{k}}{\poq{a/y}{k}}(yz)^k+\sum_{k=0}^\infty \frac{bzq^k-aG_1(zq^k)}{1-aq^k/y}\frac{\poq{b/y}{k}}{\poq{a/y}{k}}(yz)^k\\
&=\sum_{k=0}^\infty\frac{\poq{b/y}{k}}{\poq{a/y}{k}}(yz)^k+\sum_{k=0}^\infty (bzq^k-aG_1(zq^k))\frac{\poq{b/y}{k}}{\poq{a/y}{k+1}}(yz)^k.
\end{align*}
So the lemma is proved.
\qed

There also exists a finite univariate generating function of $\{B_{n,k}(a,b)\}_{k=0}^n$.
\begin{tl}\label{eee} Let $B_{n,k}(a,b)$ be defined by \eqref{18}. Then for integer $n\geq 1$, we have
\begin{align}\sum_{k=0}^nB_{n,k}(a,b)y^k=
\frac{\poq{b/y}{n-1}}{\poq{a/y}{n}}y^{n}-a\sum_{k=0}^{n-1} B_{n-k,1}(a,b)q^{(n-k)k} \frac{\poq{b/y}{k}}{\poq{a/y}{k+1}}y^k.\label{fff}
\end{align}
\end{tl}
\pf It is an immediate consequence of Lemma \ref{bgf}. 
To be precise, by Lemma \ref{bgf}, we see
\begin{align*}
G(y,z)&=\sum_{n\geq k\geq 0}B_{n,k}(a,b)y^kz^n\\
&=\sum_{n=0}^\infty\frac{\poq{b/y}{n}}{\poq{a/y}{n}}(yz)^n+\sum_{n=0}^\infty (bzq^n-aG_1(zq^n))\frac{\poq{b/y}{n}}{\poq{a/y}{n+1}}(yz)^n.
\end{align*}
By equating the coefficients of $z^n$ on both sides, it is  easy to   calculate that for  $n\geq 1$,
\begin{align*}
\sum_{k=0}^nB_{n,k}(a,b)y^k&=\frac{\poq{b/y}{n}}{\poq{a/y}{n}}y^n+b\frac{\poq{b/y}{n-1}}{\poq{a/y}{n}}(qy)^{n-1}\\
&-a\sum_{k=0}^n B_{n-k,1}(a,b)q^{(n-k)k}\frac{\poq{b/y}{k}}{\poq{a/y}{k+1}}y^k\\
&=\frac{\poq{b/y}{n-1}}{\poq{a/y}{n}}y^{n}-a\sum_{k=0}^{n-1}B_{n-k,1}(a,b)q^{(n-k)k} \frac{\poq{b/y}{k}}{\poq{a/y}{k+1}}y^k.
\end{align*} The corollary is thus proved.
\qed

\begin{remark}
Evidently, the left-hand side of \eqref{fff} is a polynomial in $y$ while the right-hand side doesn't seem that case. In fact,  the coefficients $B_{n-k,1}(a,b)$ given by \eqref{180}, i.e.,
$$
z=\sum_{n=1}^\infty B_{n,1}(a,b) z^n\frac{\poq{az}{n}}{\poq{bz}{n}},
$$
just satisfy
$$S_n(aq^k)=0\,\,\,(0\leq k\leq n-1),$$
where $S_n(y)$ is given by
\begin{align}
\frac{S_n(y)}{\prod_{k=0}^{n-1}(y-aq^k)}:=\frac{\poq{b/y}{n-1}}{\poq{a/y}{n}}y^{n}-a\sum_{k=0}^{n-1} B_{n-k,1}(a,b)q^{(n-k)k}\frac{\poq{b/y}{k}}{\poq{a/y}{k+1}}y^k.
\end{align}
This fact guarantees that the right-hand side of \eqref{fff} is  really a polynomial in $y$.
\end{remark}
Corollary \ref{eee} leads us to a general matrix inversion, which will play a very crucial role in our main result, i.e., Theorem \ref{analogue-two}.
\begin{dl}[Matrix inversion]\label{matrixinversion}
Let $A=(A_{n,k})$  be the infinitely lower-triangular matrix with the entries
\begin{align}
A_{n,k}=\boldsymbol\lbrack z^{n-k}\boldsymbol\rbrack\bigg\{\frac{\poq{az}{k}}{\poq{bz}{k}}\bigg\}\label{eq1-13-inverse}
\end{align}
and assume $A^{-1}=(B_{n,k}(a,b))$. Then
\begin{align}B_{n,k}(a,b)&=\boldsymbol\lbrack z^{n-k}\boldsymbol\rbrack\bigg\{
\frac{\poq{bz}{n-1}}{\poq{az}{n}}\bigg\}-a\sum_{i=k}^{n-1}B_{n-i,1}(a,b)q^{(n-i)i} \boldsymbol\lbrack z^{i-k}\boldsymbol\rbrack\bigg\{\frac{\poq{bz}{i}}{\poq{az}{i+1}}\bigg\}.\label{ssss}
\end{align}
\end{dl}
\pf It is clear that \eqref{ssss} is valid for $n=k=0$ or $k=0$. Thus we only need to show \eqref{ssss} for $n\geq 1$. To that end,  we first set $y\to 1/t$  in \eqref{fff} and then multiply both sides with $t^n$. All that we obtained  is
\begin{align}\sum_{k=0}^nB_{n,k}(a,b)t^{n-k}=
\frac{\poq{bt}{n-1}}{\poq{at}{n}}-a\sum_{k=0}^{n-1}B_{n-k,1}(a,b) \frac{\poq{bt}{k}}{\poq{at}{k+1}}q^{(n-k)k}t^{n-k}.\label{hhh}
\end{align}
A comparison of the coefficients of $t^{n-k}$ yields
\begin{align*}B_{n,k}(a,b)&=\boldsymbol\lbrack t^{n-k}\boldsymbol\rbrack\bigg\{
\frac{\poq{bt}{n-1}}{\poq{at}{n}}\bigg\}-a\sum_{i=0}^{n-1}B_{n-i,1}(a,b)q^{(n-i)i} \boldsymbol\lbrack t^{n-k}\boldsymbol\rbrack\bigg\{\frac{\poq{bt}{i}}{\poq{at}{i+1}}t^{n-i}\bigg\}\\
&=\boldsymbol\lbrack t^{n-k}\boldsymbol\rbrack\bigg\{
\frac{\poq{bt}{n-1}}{\poq{at}{n}}\bigg\}-a\sum_{i=k}^{n-1}B_{n-i,1}(a,b)q^{(n-i)i} \boldsymbol\lbrack t^{i-k}\boldsymbol\rbrack\bigg\{\frac{\poq{bt}{i}}{\poq{at}{i+1}}\bigg\}.
\end{align*}
Thus \eqref{ssss} is confirmed.
\qed

As byproducts of our analysis, we find two interesting properties for $\{B_{n,k}(a,b)\}_{n\geq k\geq 0}$  as follows.
\begin{xinzhi} Let $B_{n,k}(a,b)$ be given by \eqref{18}.
Then for integer $n\geq 1$ and
$t\in \mathbb{C}$, we have
\begin{align}
B_{n,k}(at,bt)&=B_{n,k}(a,b)t^{n-k},\label{maxinrong}
\\
\boldsymbol\lbrack z^{n}\boldsymbol\rbrack\bigg\{
\frac{\poq{bz}{n-1}}{\poq{az}{n}}\bigg\}&=a\sum_{i=0}^{n-1}B_{n-i,1}(a,b)q^{(n-i)i} \boldsymbol\lbrack z^{i}\boldsymbol\rbrack\bigg\{\frac{\poq{bz}{i}}{\poq{az}{i+1}}\bigg\}.\label{maxinrong-1}\end{align}
\end{xinzhi}
\pf To establish \eqref{maxinrong}, it only needs to take  $(a,b)\to (at,bt)$ in \eqref{18}. Then it follows
\begin{align*}
z^k=\sum_{n=k}^\infty\,B_{n,k}(at,bt)z^n\frac{\poq{atz}{n}}{\poq{btz}{n}}.
\end{align*}
On the other hand, on setting $z\to t\,z$ in \eqref{18}, we have
\begin{align*}
z^k=\sum_{n=k}^\infty\,B_{n,k}(a,b)t^{n-k}z^n\frac{\poq{atz}{n}}{\poq{btz}{n}}.
\end{align*}
By the uniqueness of series expansion, we obtain  \eqref{maxinrong}.
Identity \eqref{maxinrong-1} is a special case $k=0$ of \eqref{ssss}, noting that for $n\geq 1,B_{n,0}=0$.
\qed

After these preliminaries we are prepared to show
Theorem \ref{analogue-two}.

\pf The existence of \eqref{expan-two-1} is evident, because
$$
\left\{z^n\frac{\poq{az}{n}}{\poq{bz}{n}}\bigg|n=0,1,2\cdots\right\}$$
 is the base of $\mathbb{C}[[z]]$.
 Thus it only needs  to evaluate the coefficients $c_n$ in \eqref{expan-two-1}. To do this, by Theorem \ref{matrixinversion}, we have
\begin{align*}
c_n&=\sum_{k=0}^n \boldsymbol\lbrack z^{k}\boldsymbol\rbrack\left\{F(z)\right\} B_{n,k}(a,b)\\
&=\sum_{k=0}^n\boldsymbol\lbrack z^{k}\boldsymbol\rbrack\left\{F(z)\right\}\boldsymbol\lbrack z^{n-k}\boldsymbol\rbrack\bigg\{
\frac{\poq{bz}{n-1}}{\poq{az}{n}}\bigg\}\\
&-a\sum_{i=0}^{n-1}B_{n-i,1}(a,b)q^{(n-i)i} \sum_{k=0}^i\boldsymbol\lbrack z^{k}\boldsymbol\rbrack\{F(z)\}[z^{i-k}]\bigg\{\frac{\poq{bz}{i}}{\poq{az}{i+1}}\bigg\}\\
&=\boldsymbol\lbrack z^{n}\boldsymbol\rbrack\bigg\{ F(z)\frac{\poq{bz}{n-1}}{\poq{az}{n}}\bigg\}-a\sum_{i=0}^{n-1}B_{n-i,1}(a,b)q^{(n-i)i}\boldsymbol\lbrack z^{i}\boldsymbol\rbrack\bigg\{  F(z)\frac{\poq{bz}{i}}{\poq{az}{i+1}}\bigg\}.
\end{align*}
The conclusion is proved.
\qed

\begin{remark}
It is worth mentioning that in \cite{garsia-1} A.M. Garsia  and J. Remmel set up a $q$-Lagrange inversion formula, which asserts that for any formal power series $F(z)=\sum_{n=1}^{\infty}F_nz^n$ and $f(z)=\sum_{n=1}^{\infty}f_nz^n$ with
$F_1f_1\neq 0$,
\begin{eqnarray}
\sum_{n=1}^\infty f_nF(z)F(zq)\cdots F(zq^{n-1})=z\label{1}
\end{eqnarray}
holds  if and only if
\begin{eqnarray}
\sum_{n=1}^\infty F_nf(z)f(z/q)\cdots f(z/q^{n-1})=z.\label{2}
\end{eqnarray}
However, to the author's knowledge,  it is still hard to find out explicit expressions of $f_n:=B_{n,1}(a,b)$ from \eqref{1} even if $F(z)=z(1-az)/(1-bz)$.
\end{remark}
In the following, we shall examine a few specific formal expansion formulas covered by  Theorem \ref{analogue-two}. As a first consequence,  when  $a=0$ we recover  Carlitz's $q$-expansion formula {\cite[p. 206, Eq. (1.11)]{carliz}}.
\begin{tl}For any $F(z)\in \mathbb{C}[[z]]$, we have
 \begin{subequations}\label{expan-three}
\begin{equation}
F(z)=\sum_{n=0}^{\infty}\frac{c_nz^{n}}{\poq{bz}{n}},\label{expan-three-1}
\end{equation}
where
\begin{equation}
c_n=\boldsymbol\lbrack z^{n}\boldsymbol\rbrack\big\{F(z)\poq{bz}{n-1}\big\}.\label{expan-three-2}
\end{equation}
\end{subequations}
\end{tl}

We remark that Carlitz's $q$-expansion formula is a useful $q$-analogue of the  Lagrange--B\"{u}rmann inversion  formula.  The reader may consult the survey
\cite{11} of D. Stanton concerning this topic.

A second interesting consequence occurs when $b=0$.
\begin{tl} Let $B_{n,1}(a,0)$ be given by \eqref{180}. Then
\begin{subequations}\label{expan-two-11}
\begin{equation}
F(z)=\sum_{n=0}^{\infty}c_nz^{n}\poq{az}{n},\label{expan-two-22}
\end{equation}
where the coefficients
\begin{eqnarray}
c_n=\boldsymbol\lbrack z^{n}\boldsymbol\rbrack\bigg\{\frac{F(z)}{\poq{az}{n}}\bigg\}
-a\sum_{k=0}^{n-1}B_{n-k,1}(a,0)q^{(n-k)k}\boldsymbol\lbrack z^{k}\boldsymbol\rbrack\bigg\{\frac{F(z)}{\poq{az}{k+1}}\bigg\}.
\end{eqnarray}
\end{subequations}
\end{tl}

As a third consequence, the special case $b=aq$  leads us to
\begin{tl}\label{coro310} Let $F(z)=\sum_{n=0}^{\infty}a_nz^n$ and $F_k(z)=\sum_{i=0}^ka_iz^i$, being the $k$-truncated series of $F(z)$. Suppose that
 \begin{subequations}\label{expan-four}
\begin{equation}
\frac{F(z)}{1-az}=\sum_{n=0}^{\infty}\frac{c_nz^{n}}{1-azq^n}.\label{expan-four-1-0}
\end{equation}
Then $c_0=a_0$ and $n\geq1,$
\begin{equation}
c_n=\sum_{k=0}^{n-1}g_{n-k}(q)q^{(n-k)k}\boldsymbol\lbrack z^{n}\boldsymbol\rbrack\bigg\{  \frac{F(z)-F_k(z)}{1-az}\bigg\},\label{expan-four-2-0}
\end{equation}
\end{subequations}
where  $g_n(q)$ are  polynomials in $q$ given recursively by
\begin{align*}g_n(q)=1
-\sum_{i=1}^{n-1}g_{n-i}(q)q^{(n-i)i}.
\end{align*}
\end{tl}
\pf In such case, we first solve  the recurrence relation \eqref{ssss} with $k=1$ for $B_{n,1}(a,aq)$, viz.
\begin{align*}B_{n,1}(a,aq)=a^{n-1}
-\sum_{i=1}^{n-1}B_{n-i,1}(a,aq)q^{(n-i)i} a^{i}.
\end{align*}
The solution is recursively given  by
\begin{align}
\left\{
  \begin{array}{ll}
 &B_{n,1}(a,aq)=g_n(q)a^{n-1},\\
&\\
&g_n(q)=\displaystyle1-\sum_{i=1}^{n-1}g_{n-i}(q)q^{(n-i)i}.\label{crucial}
  \end{array}
\right.
\end{align}
By virtue of \eqref{crucial}, we are now able to calculate  $c_n$. To do this, by Theorem \ref{analogue-two} we have
\begin{align*}
c_n&=\boldsymbol\lbrack z^{n}\boldsymbol\rbrack\bigg\{  \frac{F(z)}{1-az}\bigg\}-\sum_{k=0}^{n-1}g_{n-k}(q)q^{(n-k)k}a^{n-k}\boldsymbol\lbrack z^{k}\boldsymbol\rbrack\bigg\{\frac{F(z)}{1-az}\bigg\}\\
&=\boldsymbol\lbrack z^{n}\boldsymbol\rbrack\bigg\{  \frac{F(z)}{1-az}\bigg\}-\sum_{k=0}^{n-1}g_{n-k}(q)q^{(n-k)k}\boldsymbol\lbrack z^{n}\boldsymbol\rbrack\bigg\{\frac{F_k(z)}{1-az}\bigg\}\\
&=\boldsymbol\lbrack z^{n}\boldsymbol\rbrack\bigg\{\sum_{k=0}^{n-1}g_{n-k}(q)q^{(n-k)k}\frac{F(z)-F_k(z)}{1-az}\bigg\}.
\end{align*}
In  the penultimate equality we have used the fact that
\begin{align*}
a^{n-k}\boldsymbol\lbrack z^{k}\boldsymbol\rbrack\bigg\{\frac{F(z)}{1-az}\bigg\}=\boldsymbol\lbrack z^{n}\boldsymbol\rbrack\bigg\{\frac{F_k(z)}{1-az}\bigg\}
\end{align*}
and in the last equality, we have invoked \eqref{crucial} again.
The conclusion is proved.
\qed

It is also of interest to note that if $F(z)$ is a polynomial of degree $m+1$, say
$$
F(z)=(1-az)\prod_{n=1}^{m}(1-t_nz),
$$
and  $F_k(z)=F(z)$ for $k\geq m+1$, then Corollary \ref{coro310} reduces to
\begin{tl} With the same notation as Corollary \ref{coro310}. Then we have
 \begin{subequations}
\begin{equation}
\prod_{n=1}^{m}(1-t_nz)=\sum_{n=0}^\infty \frac{c_nz^{n}}{1-azq^n},\label{exexpan-two-1-0}
\end{equation}
where $c_0=1$ and $n\geq 1$,
\begin{equation}
c_n=\sum_{k=0}^{\min\{m,n-1\}} g_{n-k}(q)q^{(n-k)k}\boldsymbol\lbrack z^{n}\boldsymbol\rbrack\bigg\{\prod_{i=1}^{m}(1-t_iz)-\frac{F_k(z)}{1-az}\bigg\}.\label{mistake-wangjin}
\end{equation}
\end{subequations}
\end{tl}
\section{Applications to $q$-series theory}
Unlike the preceding section, we now focus our attention on applications of Theorem \ref{analogue-two}  to  the $q$-series theory. In this sense, we assume that all results are subject to appropriate convergent conditions of rigorous analytic theory, unless otherwise stated.

Let us begin with the proof of  Theorem \ref{tlidentity}.

\pf We only need to make use of Theorem \ref{analogue-two} as well as the $q$-binomial theorem \cite[(II.3)]{10} to get
\begin{align*}
\frac{\poq{az}{\infty}}{\poq{bz}{\infty}}G(z)=\sum_{n=0}^{\infty}c_nz^{n}\frac{\poq{az}{n}}{\poq{bz}{n}}
=S_1-aS_2,
\end{align*}
where
\begin{align*}
 S_1&:=\sum_{n=0}^{\infty}\sum_{i=0}^{n}\frac{\poq{aq/b}{i}}{\poq{q}{i}}b^iq^{(n-1)i}t_{n-i}z^{n}\frac{\poq{az}{n}}{\poq{bz}{n}},\\
 S_2&:=\sum_{n=0}^{\infty}\sum_{k=0}^{n-1}B_{n-k,1}(a,b)q^{(n-k)k}\sum_{i=0}^k
\frac{\poq{aq/b}{i}}{\poq{q}{i}}(bq^k)^it_{k-i}z^{n}\frac{\poq{az}{n}}{\poq{bz}{n}}.
\end{align*}
After a mere series rearrangement, we get
\begin{align*}
 S_1&=\sum_{i=0}^{\infty}\frac{\poq{aq/b}{i}}{\poq{q}{i}}b^iq^{i^2-i}z^{i}\frac{\poq{az}{i}}{\poq{bz}{i}}\sum_{n= 0}^{\infty}q^{ni}t_{n}z^{n}\frac{\poq{azq^i}{n}}{\poq{bzq^i}{n}}\\
 &=\sum_{i= 0}^{\infty}\frac{\poq{aq/b,az}{i}}{\poq{q,bz}{i}}b^iq^{i^2-i}z^{i}\widetilde{G}(zq^i;a,b).
\end{align*}
Hereafter, as given by \eqref{dy},
\[\widetilde{G}(z;a,b)=\sum_{n=0}^{\infty}t_{n}z^{n}\frac{\poq{az}{n}}{\poq{bz}{n}}.\]
In a similar way, it is easily found that
\begin{align*}
  S_2&=\sum_{i=0}^{\infty}\frac{\poq{aq/b}{i}}{\poq{q}{i}}b^i
\sum_{k=i}^{\infty}t_{k-i}z^k\frac{\poq{az}{k}}{\poq{bz}{k}}q^{ki}\Delta_{k},
\end{align*}
where
\[
\Delta_{k}:=\sum_{n= k+1}^{\infty}B_{n-k,1}(a,b)(zq^{k})^{n-k}\frac{\poq{azq^k}{n-k}}{\poq{bzq^k}{n-k}}=\sum_{n= 1}^{\infty}B_{n,1}(a,b)(zq^{k})^{n}\frac{\poq{azq^k}{n}}{\poq{bzq^k}{n}}=
zq^k.
\]
The last equality is based on \eqref{180}.
Therefore,
\begin{align*}
  S_2&=z\sum_{i=0}^{\infty}\frac{\poq{aq/b,az}{i}}{\poq{q,bz}{i}}(bzq^{i+1})^i
\sum_{k=0}^{\infty}t_{k}\frac{\poq{azq^i}{k}}{\poq{bzq^i}{k}}(zq^{i+1})^{k}\\
&=z\sum_{i=0}^{\infty}\frac{\poq{aq/b,az}{i}}{\poq{q,bz}{i}}(bzq^{i+1})^i\widetilde{G}(zq^{i+1};a/q,b/q).
\end{align*}
Finally, we achieve
\begin{align*}
\frac{\poq{az}{\infty}}{\poq{bz}{\infty}}G(z)&=
\sum_{n=0}^{\infty}\frac{\poq{aq/b,az}{n}}{\poq{q,bz}{n}}(bz)^nq^{n(n-1)}\big(\widetilde{G}(zq^n;a,b)-
azq^{2n}\widetilde{G}(zq^{n+1};a/q,b/q)\big).
\end{align*}
This gives the complete proof of the theorem.
\qed

With regard to applications of  Theorem  \ref{tlidentity} to $q$-series, it is necessary to set up
\begin{tl}\label{tlnew}For integer $r\geq 0$ and $|cz|<1$, it holds
\begin{align}
&\frac{\poq{azq}{\infty}}{\poq{bz}{\infty}}{_{r+1}\phi_r}\bigg[\genfrac{}{}{0pt}{}{A_1,A_2,\ldots,A_{r+1}}{B_1,B_2,\ldots,B_{r}};q,cz\bigg]\nonumber\\
&=\lim_{x,y\to\infty}\sum_{n=0}^{\infty}\frac{\poq{az,(az)^{1/2}q,-(az)^{1/2}q,aq/b,x,y}{n}}{\poq{q,(az)^{1/2},-(az)^{1/2},bz,azq/x,azq/y}{n}}\bigg(\frac{bz}{xy}\bigg)^n\nonumber\\
&\qquad\quad\times{_{r+3}\phi_{r+2}}\bigg[\genfrac{}{}{0pt}{}{azq^n,azq^{2n+1},A_1,A_2,\ldots,A_{r+1}}{bzq^n,azq^{2n},B_1,B_2,\ldots,B_{r}};q,czq^n\bigg].\label{expan-two-5lz}
\end{align}
\end{tl}
\pf It suffices to set in Theorem \ref{tlidentity}
$$G(z)={_{r+1}\phi_r}\bigg[\genfrac{}{}{0pt}{}{A_1,A_2,\ldots,A_{r+1}}{B_1,B_2,\ldots,B_{r}};q,cz\bigg],$$ which means
 $$t_k=\frac{\poq{A_1,A_2,\ldots,A_{r+1}}{k}}{\poq{q,B_1,B_2,\ldots,B_{r}}{k}}c^k.$$
In the sequel, it is routine to compute
\begin{align*}
H_n(z;a,b)&:=\frac{\widetilde{G}(zq^n;a,b)-
azq^{2n}\widetilde{G}(zq^{n+1};a/q,b/q)}{1-azq^{2n}}\\
&=\sum_{k=0}^\infty\frac{\poq{azq^n,azq^{2n+1}}{k}}{\poq{bzq^n,azq^{2n}}{k}}(zq^n)^k t_k\\
&={_{r+3}\phi_{r+2}}\bigg[\genfrac{}{}{0pt}{}{azq^n,azq^{2n+1},A_1,A_2,\ldots,A_{r+1}}{bzq^n,azq^{2n},B_1,B_2,\ldots,B_{r}};q,czq^n\bigg].
\end{align*}
This reduces  \eqref{expan-two-lzlzlz} of Theorem \ref{tlidentity} to
\begin{align}
&\frac{\poq{azq}{\infty}}{\poq{bz}{\infty}}{_{r+1}\phi_r}\bigg[\genfrac{}{}{0pt}{}{A_1,A_2,\ldots,A_{r+1}}{B_1,B_2,\ldots,B_{r}};q,cz\bigg]\nonumber\\
&\qquad=
\sum_{n=0}^{\infty}\frac{\poq{aq/b,az}{n}}{\poq{q,bz}{n}}(bz)^nq^{n(n-1)}\frac{1-azq^{2n}}{1-az}H_{n}(z;a,b).\label{xuyao}
\end{align}
Finally, using the basic relations
\begin{align*}
\frac{1-azq^{2n}}{1-az}
=\frac{\poq{(az)^{1/2}q,-(az)^{1/2}q}{n}}{\poq{(az)^{1/2},-(az)^{1/2}}{n}}
\end{align*}
and
\begin{align*}
\lim_{x,y\to\infty}
\frac{\poq{x,y}{n}}{\poq{azq/x,azq/y}{n}}\bigg(\frac{1}{xy}\bigg)^n
=q^{n(n-1)},
\end{align*}
we derive \eqref{expan-two-5lz} from \eqref{xuyao} directly.
\qed

The following are two special instances of Theorem \ref{tlidentity}.
\begin{lz}The following transformation formulas are valid.
\begin{align}
\frac{\poq{az}{\infty}}{\poq{bz}{\infty}}&=
\sum_{n=0}^{\infty}\frac{\poq{aq/b,az}{n}}{\poq{q,bz}{n}}(bz)^nq^{n(n-1)}(1-azq^{2n}),\label{expan-two-lzlz}\\
\frac{\poq{azq,ABz}{\infty}}{\poq{bz,Bz}{\infty}}&=
\sum_{n=0}^{\infty}\frac{\poq{(az)^{1/2}q,-(az)^{1/2}q,aq/b,az}{n}}{\poq{q,(az)^{1/2},-(az)^{1/2},bz}{n}}
(bz)^nq^{n(n-1)}\nonumber\\
&\qquad\times{_3\phi_2}\bigg[\genfrac{}{}{0pt}{}{azq^n,azq^{2n+1},A}{bzq^n,azq^{2n}};q,Bzq^n\bigg].\label{expan-two-4lz}
\end{align}
\end{lz}
\pf Identity  \eqref{expan-two-lzlz}  comes from $G(z)=1$ in Theorem \ref{tlidentity} and  \eqref{expan-two-4lz} does by taking $G(z)=\poq{ABz}{\infty}/\poq{Bz}{\infty}$, i.e., $t_k=\poq{A}{k}/\poq{q}{k}B^k$
in Theorem \ref{tlidentity} or $r=0$ in Corollary \ref{tlnew}.
\qed

The next conclusion shows how  Theorem \ref{tlidentity} can  be applied to known transformation formulas for finding new results.
\begin{tl} For $|z|<1$, we have
\begin{align} {}_{2}\phi _{1}\left[\begin{matrix}A,B\\ C\end{matrix}
; q, z\right]=&\sum_{n=0}^{\infty}
\frac{\poq{ABq/C,ABz/C}{n}}
{\poq{q,z}{n}} z^n q^{n(n-1)}\left(1-\frac{ABzq^{2n}}{C}\right)\nonumber\\
&\times{}_{4}\phi _{3}\left[\begin{matrix}ABzq^n/C,ABzq^{2n+1}/C,C/A,C/B\\ C,zq^n,ABzq^{2n}/C\end{matrix}
; q, \frac{ABzq^n}{C}\right].
\end{align}
\end{tl}
\pf  Performing as above, we choose in Theorem \ref{tlidentity}
\begin{align*}
G(z)={}_{2}\phi _{1}\left[\begin{matrix}C/A,C/B\\ C\end{matrix}
; q, \frac{ABz}{C}\right],\end{align*}
which corresponds to
\begin{align*}
t_k=\frac{\poq{C/A,C/B}{k}}{\poq{q,C}{k}}\left(\frac{AB}{C}\right)^k.
\end{align*}In this case, it is clear that
\begin{align*}
H_n(z;a,b):&=\frac{\widetilde{G}(zq^n;a,b)-
azq^{2n}\widetilde{G}(zq^{n+1};a/q,b/q)}{1-azq^{2n}}\\
&={}_{4}\phi _{3}\left[\begin{matrix}azq^n,azq^{2n+1},C/A,C/B\\ C,bzq^n,azq^{2n}\end{matrix}
; q, \frac{ABzq^n}{C}\right].
\end{align*}
As a result, from Theorem \ref{tlidentity} it follows
\begin{align*}
\frac{\poq{az}{\infty}}{\poq{bz}{\infty}}{}_{2}\phi _{1}\left[\begin{matrix}C/A,C/B\\ C\end{matrix}
; q, \frac{ABz}{C}\right]&=
\sum_{n=0}^{\infty}\frac{\poq{aq/b,az}{n}}{\poq{q,bz}{n}}(bz)^nq^{n(n-1)}(1-azq^{2n})H_n(z;a,b).
\end{align*}
In this form, taking  $a=AB/C$ and $b=1$, we  obtain
\begin{align}
&\frac{\poq{ABz/C}{\infty}}{\poq{z}{\infty}}{}_{2}\phi _{1}\left[\begin{matrix}C/A,C/B\\ C\end{matrix}
; q, \frac{ABz}{C}\right]\nonumber\\
&=\sum_{n=0}^{\infty}\frac{\poq{ABq/C,ABz/C}{n}}{\poq{q,z}{n}}z^nq^{n(n-1)}(1-ABzq^{2n}/C)H_n(z;AB/C,1).\label{need}
\end{align}
By combining \eqref{need} with Heine's third transformation \cite[(III.3)]{10}
\begin{align*}
{}_{2}\phi _{1}\left[\begin{matrix}A,B\\ C\end{matrix}
; q, z\right]=&\frac{\poq{ABz/C}{\infty}}{\poq{z}{\infty}}{}_{2}\phi _{1}\left[\begin{matrix}C/A,C/B\\ C\end{matrix}
; q, \frac{ABz}{C}\right],
\end{align*}
then  we reformulate  \eqref{need} in standard notation of $q$-series  as
\begin{align*} {}_{2}\phi _{1}\left[\begin{matrix}A,B\\ C\end{matrix}
; q, z\right]=&\sum_{n=0}^{\infty}
\frac{\poq{ABq/C,ABz/C}{n}}
{\poq{q,z}{n}} z^n q^{n(n-1)}\left(1-\frac{ABzq^{2n}}{C}\right)\nonumber\\
&\times{}_{4}\phi _{3}\left[\begin{matrix}ABzq^n/C,ABzq^{2n+1}/C,C/A,C/B\\ C,zq^n,ABzq^{2n}/C\end{matrix}
; q, \frac{ABzq^n}{C}\right].
\end{align*}
The conclusion is proved.
\qed

Perhaps, the most interesting case  is the following partial theta function identity. It can be derived from
  Theorem \ref{tlidentity} with the help of two Coogan-Ono type identities \eqref{id11} and \eqref{id1177}.
\begin{tl}[Partial theta function identity]Let $\theta(z;q)$ be the partial theta function given by
\[\sum_{n=0}^\infty(-1)^nq^{n(n-1)/2}z^{n}.\] Then
\begin{align}
&\frac{\poq{zq}{\infty}}{\poq{-zq}{\infty}}+\sum_{n=0}^{\infty}\frac{\poq{-1,z}{n}}{\poq{q,-zq}{n}}(-z)^nq^{n^2+n}
\label{expan-two-lzlzlz-new}
\\
&\qquad=\sum_{n=0}^{\infty}\frac{\poq{-1,z}{n}}{\poq{q,-zq}{n}}\big(1+q^n+zq^{n}-zq^{2n}\big)(-z)^nq^{n^2}\theta(z^2q^{2n+1};q^2).\nonumber
\end{align}
\end{tl}
\pf Recall that Lemma \ref{ma-id11} gives
\begin{align}
\sum_{k=0}^\infty z^{k}\frac{\poq{z}{k}}{\poq{-zq}{k}}=(1+z)\sum_{k=0}^\infty(-1)^kz^{2k}q^{k^2}.\label{idnewadded-1}
\end{align}
Lemma \ref{ma-id1177} can be restated as
\begin{align}
\sum_{k=0}^\infty\,z^k\frac{\poq{z}{k}}{\poq{-zq}{k}}(1-zq^k)=1+
2\sum_{k=1}^{\infty}(-1)^kz^{2k}q^{k^2}.\label{idnewadded-0}
\end{align}
Subtracting  \eqref{idnewadded-0} from \eqref{idnewadded-1}, we obtain
\begin{align}
z\sum_{k=0}^\infty\,(qz)^k\frac{\poq{z}{k}}{\poq{-zq}{k}}&=z+(z-1)\sum_{k=1}^{\infty}(-1)^kz^{2k}q^{k^2},\nonumber\\
\mbox{i.e.,}\,\,\,\sum_{k=0}^\infty\,(qz)^k\frac{\poq{z}{k}}{\poq{-zq}{k}}&=\sum_{k=0}^{\infty}(-1)^kz^{2k}q^{k^2}-\sum_{k=1}^{\infty}(-1)^kz^{2k-1}q^{k^2}.
\label{idnewadded-2}
\end{align}
Using \eqref{idnewadded-1} and \eqref{idnewadded-2},  as well as referring to \eqref{dy} with $t_k=1$, we thus obtain
\begin{align}
\widetilde{G}(z;1,-q)&=(1+z)\sum_{k=0}^\infty(-1)^kz^{2k}q^{k^2},\\
\widetilde{G}(zq;1/q,-1)&=\sum_{k=0}^{\infty}(-1)^kz^{2k}q^{k^2}-\sum_{k=1}^{\infty}(-1)^kz^{2k-1}q^{k^2}.
\end{align}
Thus it is easy to check that
\begin{align*}\displaystyle
&\widetilde{G}(zq^n;1,-q)-
zq^{2n}\widetilde{G}(zq^{n+1};1/q,-1)\\
&=(1+zq^n)\sum_{k=0}^\infty(-1)^kz^{2k}q^{k^2+2kn}-zq^{2n}\sum_{k=0}^{\infty}(-1)^kz^{2k}q^{k^2+2kn}+q^n\sum_{k=1}^{\infty}(-1)^kz^{2k}q^{k^2+2kn}\\
&=-q^n+(1+q^n+zq^n-zq^{2n}
)\sum_{k=0}^{\infty}(-1)^kz^{2k}q^{k^2+2kn}.
\end{align*}
Note that the summation on the right-hand side can be recast in terms of $\theta(z;q)$. We thus obtain
\begin{align*}\displaystyle
&\widetilde{G}(zq^n;1,-q)-
zq^{2n}\widetilde{G}(zq^{n+1};1/q,-1)
=-q^n+(1+q^n+zq^n-zq^{2n}
)\theta(z^2q^{2n+1};q^2).
\end{align*}
This reduces the whole equation \eqref{expan-two-lzlzlz} to
\begin{align*}
\frac{\poq{zq}{\infty}}{\poq{-zq}{\infty}}&+\sum_{n=0}^{\infty}\frac{\poq{-1,z}{n}}{\poq{q,-zq}{n}}(-z)^nq^{n^2+n}
\\
&=\sum_{n=0}^{\infty}\frac{\poq{-1,z}{n}}{\poq{q,-zq}{n}}(-z)^nq^{n^2}\big(1+q^n+zq^{n}-zq^{2n}\big)\theta(z^2q^{2n+1};q^2).
\end{align*}
Thus Identity \eqref{expan-two-lzlzlz-new} is proved.
\qed

In the case that $z=q^{-m},m\geq 1$, \eqref{expan-two-lzlzlz-new} reduces to  a finite summation of $\theta(z;q).$
\begin{lz}For $m\geq 1$, we have
\begin{align}
&\sum_{n=0}^{m}\frac{\poq{-1}{n}}{\poq{-q^{1-m}}{n}}
\bigg[\genfrac{}{}{0pt}{}{m}{n}\bigg]_qq^{3n^2/2+n/2-2nm}\label{qqq}\\
&=\sum_{n=0}^{m}\frac{\poq{-1}{n}}{\poq{-q^{1-m}}{n}}
\bigg[\genfrac{}{}{0pt}{}{m}{n}\bigg]_qq^{3n^2/2-n/2-2nm}
\big(1+q^n+q^{n-m}-q^{2n-m}\big)\theta(q^{2n-2m+1};q^2),\nonumber
\end{align}
where $\bigg[\genfrac{}{}{0pt}{}{m}{n}\bigg]_q$ is the usual $q$-binomial coefficient.
\end{lz}
\pf It suffices to take $z=q^{-m}$ in \eqref{expan-two-lzlzlz-new} and simplify the obtained by using the facts that for integer $m\geq 1$, $\poq{q^{1-m}}{\infty}=0$ and
\begin{align*}
\frac{\poq{q^{-m}}{n}}{\poq{q}{n}}=
\bigg[\genfrac{}{}{0pt}{}{m}{n}\bigg]_q(-1)^nq^{n(n-1)/2-mn}.
\end{align*}
\qed

It would be natural to expect that  Theorem \ref{analogue-two} can be applied to bilateral $q$-series. The reader is referred to  \cite[Eq. (5.1.2)]{10} or \eqref{guiding} for the definition of bilateral $q$-series. As an interesting example,  we now set up  a coefficient identity of the famous Ramanujan ${}_1\psi_1$ summation formula \cite[(II.29)]{10}.

\begin{tl}Let $B_{n,1}(a,b)$ be given by \eqref{180}. For  $|b/a|<|z|<1,$ and integer $n\geq 0$, it holds
\begin{align}
&\boldsymbol\lbrack z^{n}\boldsymbol\rbrack\bigg\{\frac{(aqz^2,1/az^2;q)_{\infty}}
{(z,b/az,bq^nz,1/az;q)_{\infty}\poq{aqz}{n}}\bigg\}\label{313}=\frac{1}{(q,b/a;q)_{\infty}}\\
&\quad+aq\sum_{k=0}^{n-1}B_{n-k,1}(aq,bq)q^{(n-k)k}\boldsymbol\lbrack z^{k}\boldsymbol\rbrack\bigg\{\frac{(aqz^2,1/az^2;q)_{\infty}}
{(z,b/az,bq^{k+1}z,1/az;q)_{\infty}\poq{aqz}{k+1}}\bigg\}.\nonumber
\end{align}
\end{tl}
\pf
Observe that  Ramanujan's $\,_1\psi_1$ sum   states that for  $|b/a|<|z|<1,$
\begin{eqnarray}\sum_{k=-\infty}^{\infty}\frac{(a;q)_k}{(b;q)_k}z^k=
\frac{(az,q/az,q,b/a;q)_{\infty}}
{(z,b/az,b,q/a;q)_{\infty}}.
 \end{eqnarray}
Set $(a,b)\to (aqz,bqz)$. Then we arrive at
\begin{eqnarray}\sum_{k=-\infty}^{\infty}\frac{(aqz;q)_k}{(bqz;q)_k}z^k=
\frac{(aqz^2,1/az^2,q,b/a;q)_{\infty}}
{(z,b/az,bqz,1/az;q)_{\infty}},
 \end{eqnarray}
 which can be reformulated  in the form
  \begin{align}
  f(z)+g(1/z)=F(z), \label{bilateral}
  \end{align}
 where we define
 \begin{align*}
f(z)&:=\sum_{k=0}^{\infty}\frac{(aqz;q)_k}{(bqz;q)_k}z^k,\,\,g(z):=
 \sum_{k=1}^{\infty}\frac{(z/b;q)_k}{(z/a;q)_k}\bigg(\frac{bz}{a}\bigg)^k,\\
F(z)&:=\frac{(aqz^2,1/az^2,q,b/a;q)_{\infty}}
{(z,b/az,bqz,1/az;q)_{\infty}}.
 \end{align*}
 Now we apply the expansion formula in Theorem \ref{analogue-two} to $f(z)$. It follows from \eqref{coeformula} that for $n\geq 0$,
 \begin{align}
1=\boldsymbol\lbrack z^{n}\boldsymbol\rbrack\bigg\{  f(z)\frac{\poq{bqz}{n-1}}{\poq{aqz}{n}}\bigg\}-aq\sum_{k=0}^{n-1}B_{n-k,1}(aq,bq)q^{(n-k)k}\boldsymbol\lbrack z^{k}\boldsymbol\rbrack\bigg\{f(z)\frac{\poq{bqz}{k}}{\poq{aqz}{k+1}}\bigg\}.\label{coeformula1890}
\end{align}
Next, observe that
\begin{align*}
\boldsymbol\lbrack z^{n}\boldsymbol\rbrack\bigg\{f(z)\frac{\poq{bqz}{n-1}}{\poq{aqz}{n}}\bigg\}=\sum_{k=0}^n\boldsymbol\lbrack z^{k}\boldsymbol\rbrack\{f(z)\}\times\boldsymbol\lbrack z^{n-k}\boldsymbol\rbrack\bigg\{  \frac{\poq{bqz}{n-1}}{\poq{aqz}{n}}\bigg\},
 \end{align*}
 while for $k\geq 0$, due to \eqref{bilateral}, it holds
$$\boldsymbol\lbrack z^{k}\boldsymbol\rbrack\{f(z)\}=\boldsymbol\lbrack z^{k}\boldsymbol\rbrack\{F(z)\}.$$
We immediately obtain
\begin{align*}
\boldsymbol\lbrack z^{k}\boldsymbol\rbrack\bigg\{  f(z)\frac{\poq{bqz}{k-1}}{\poq{aqz}{k}}\bigg\}=\boldsymbol\lbrack z^{k}\boldsymbol\rbrack\bigg\{F(z)\frac{\poq{bqz}{k-1}}{\poq{aqz}{k}}\bigg\}.
\end{align*}
This simplifies \eqref{coeformula1890} to
\begin{align*}
1&=\boldsymbol\lbrack z^{n}\boldsymbol\rbrack\bigg\{  F(z)\frac{\poq{bqz}{n-1}}{\poq{aqz}{n}}\bigg\}-aq\sum_{k=0}^{n-1}B_{n-k,1}(aq,bq)q^{(n-k)k}\boldsymbol\lbrack z^{k}\boldsymbol\rbrack\bigg\{F(z)\frac{\poq{bqz}{k}}{\poq{aqz}{k+1}}\bigg\}\\
&=\boldsymbol\lbrack z^{n}\boldsymbol\rbrack\bigg\{\frac{(aqz^2,1/az^2,q,b/a;q)_{\infty}}
{(z,b/az,bq^nz,1/az;q)_{\infty}\poq{aqz}{n}}\bigg\}\\
&-aq\sum_{k=0}^{n-1}B_{n-k,1}(aq,bq)q^{(n-k)k}\boldsymbol\lbrack z^{k}\boldsymbol\rbrack\bigg\{\frac{(aqz^2,1/az^2,q,b/a;q)_{\infty}}
{(z,b/az,bq^{k+1}z,1/az;q)_{\infty}\poq{aqz}{k+1}}\bigg\}.
\end{align*}
It turns out to be \eqref{313}.
\qed

We conclude our paper with a coefficient identity of the Coogan-Ono identity \eqref{id11}  which can be  easily derived by using  \eqref{18}.
\begin{tl}
Let $B_{n,k}(a,b)$ be given by \eqref{18}. Then for any integer $n\geq 0$, we have
\begin{align}
\sum_{k=0}^{\lfloor n/2\rfloor}B_{n,2k}(1,-q)(-1)^kq^{k^2}+\sum_{k=0}^{\lfloor (n-1)/2\rfloor}B_{n,2k+1}(1,-q)(-1)^kq^{k^2}=1,
\end{align}
where $\lfloor x\rfloor$ denotes the usual floor function.
\end{tl}
\pf It suffices to take $a=1,b=-q$ in Theorem \ref{analogue-two} and
$$F(z)=(1+z)\sum_{n=0}^{\infty}(-1)^nz^{2n}q^{n^2}:=\sum_{n=0}^{\infty}a_nz^n.$$
So we are back with the series expansion
\begin{align*}
F(z)=\sum_{n=0}^\infty\,z^n\frac{\poq{z}{n}}{\poq{-zq}{n}}.
\end{align*}
Thus, by  \eqref{18} instead of \eqref{coeformula}, we obtain
\begin{align*}
1=\sum_{k=0}^nB_{n,k}(1,-q)a_{k}&=\sum_{2k=0}^nB_{n,2k}(1,-q)(-1)^kq^{k^2}\\
&+\sum_{2k+1=1}^nB_{n,2k+1}(1,-q)(-1)^kq^{k^2}.
\end{align*}
\qed

 \section*{Acknowledgements}
This  work was supported by the National Natural Science Foundation of
China [Grant No.  11471237].
\bibliographystyle{amsplain}

\end{document}